\documentclass[coll]{imsart}

\usepackage{amssymb,amsmath,amsfonts,amsthm}
\providecommand{\url}[1]{}

\numberwithin{equation}{section}

\newtheorem{theorem}{Theorem}

\newtheorem{corollary}[theorem]{Corollary}

\newtheorem{definition}[theorem]{Definition}

\newtheorem{lemma}[theorem]{Lemma}

\newtheorem{remark}[theorem]{Remark}

\pubyear{2009}
\volume{0}
\volumetitle{Volume Title}
\arxiv{math.PR/0000000}

\begin{document}

\begin{frontmatter}

\title{Bernstein inequality and moderate deviations under strong
mixing conditions
}
\runtitle{Bernstein inequality under strong mixing}

\begin{aug}
\author{Florence Merlev\`ede\thanksref{a}\ead[label=e1]{florence.merlevede@univ-mlv.fr}}
\author{Magda Peligrad\thanksref{b,t1}\ead[label=e2]{magda.peligrad@uc.edu}}
\and
\author{Emmanuel Rio\thanksref{c,t2}\ead[label=e3]{rio@math.uvsq.fr}}
  

  
  \thankstext{t1}{Supported in part by a Charles Phelps Taft Memorial Fund grant and NSA\
grant, H98230-07-1-0016 and H98230-09-1-0005.}
\thankstext{t2}{Supported in part by Centre INRIA Bordeaux Sud-Ouest $\&$ Institut de
Math\'ematiques de Bordeaux}

  \runauthor{F. Merlev\`ede, M. Peligrad and E. Rio}

  \affiliation{ Universit\'e Paris Est, University of Cincinnati, Universit\'e de Versailles}

  \address[a]{Universit\'e Paris-Est-Marne-la-Vall\'{e}e, LAMA and
CNRS UMR 8050, \printead{e1}}

\address[b]{Department of Mathematical Sciences, University of Cincinnati, \printead{e2}}

\address[c]{ Universit\'e de Versailles, Laboratoire de math\'ematiques, UMR 8100
CNRS, \printead{e3}}

\contributor{Merlev\`ede, F.}{Universit\'e Paris Est-Marne la vall\'ee}
\contributor{Peligrad, M.}{University of Cincinnati}
\contributor{Rio, E.}{Universit\'e de Versailles}
\end{aug}

\begin{abstract}
In this paper we obtain a Bernstein type inequality for a class of weakly
dependent and bounded random variables. The proofs lead to a moderate deviations
principle for sums of bounded random variables with exponential decay of the
strong mixing coefficients that complements the large deviation result
obtained by Bryc and Dembo (1998) under superexponential mixing rates.
\end{abstract}




\end{frontmatter}

\section{Introduction}

This paper has double scope. First we obtain a Bernstein's type bound on the
tail probabilities of the partial sums $S_{n}$ of a sequence of dependent and bounded 
random variables $(X_{k},k\geq 1).$ Then we use the developed techniques to
study the moderate deviations principle.

We recall the definition of strongly mixing sequences, introduced by
Rosenblatt (1956): For any two $\sigma$-algebras $\mathcal{A}$ and 
$\mathcal{B}$, we define the $\alpha $-mixing coefficient by 

\begin{equation*}
\alpha (\mathcal{A},\mathcal{B})=\sup_{A\in \mathcal{A},B\in \mathcal{B}}|%
\mathbb{P}(A\cap B)-\mathbb{P}(A)\mathbb{P}(B)|\text{ }.
\end{equation*}

Let $(X_{k},k\geq 1)$ be a sequence of real-valued random variables defined on $%
\left( \Omega ,\mathcal{A},\mathbb{P}\right) $. This sequence will be called
strongly mixing if 
\begin{equation}
\alpha (n):=\sup_{k\geq 1}\alpha \left( \mathcal{M}_{k},\mathcal{G}%
_{k+n}\right) \rightarrow 0\mbox{ as }n\rightarrow \infty \,,
\label{defalpha}
\end{equation}%
where $\mathcal{M}_{j}:=\sigma (X_{i},i\leq j)$ and $\mathcal{G}_{j}:=\sigma
(X_{i},i\geq j)$ for $j\geq 1$.

Alternatively (see Bradley, 2007, Theorem 4.4) 
\begin{equation}
4\alpha (n):=\sup\{ \mathrm{Cov}(f,g)/||f||_{\infty } ||g||_{\infty };\text{ 
}f\in \mathbb{L}_{\infty }(\mathcal{M}_{k}),\text{ }g\in \mathbb{L}_{\infty
}(\mathcal{G}_{k+n}) \} \, .  \label{defalpha2}
\end{equation}%
Establishing exponential inequalities for strongly mixing sequences is a
very challenging problem. Some steps in this direction are results by Rio
(2000, Theorem 6.1), who obtained a  Fuk-Nagaev type inequality, by
Dedecker and Prieur (2004) who extended Theorem 6.1 in Rio (2000) using coupling
coefficients, by  Doukhan and Neumann (2005) who used combinatorics
techniques. In a recent paper Merlev\`{e}de, Peligrad and Rio (2009) get
exponential bounds for subexponential mixing rates, when the variables are
not necessarily bounded, obtaining the same order of magnitude as in the
independent case. More precisely they show that, if $\alpha (n)\leq \exp
(-cn^{\gamma _{1}})$ and $\sup_{i>0}\mathbb{P}(|X_{i}|>t)\leq \exp
(1-t^{\gamma _{2}})$ with $\gamma _{1}>0$ and $\gamma _{2}>0$, such that $%
(1/\gamma _{1})+(1/\gamma _{2})=1/\gamma >1$, then there are positive
constants $C$, $C_{1},$ $C_{2}$ and $\eta $ depending only on $c$, $\gamma
_{1}$ and $\gamma _{2}$, such that for all $n\geq 4$ and $\lambda \geq
C(\log n)^{\eta }$ 
\begin{equation*}
\mathbb{P}(|S_{n}|\geq \lambda )\leq (n +1) \exp (-\lambda ^{\gamma}/C_1 ) + \exp  (-\lambda ^{2}/nC_2 )\,.
\end{equation*}%
Here $S_{n}=\sum\nolimits_{k=1}^{n}X_{k}.$ The case not covered by that
paper is the case of exponential mixing rates and bounded variables, that is 
$\gamma _{1}=1$ and $\gamma _{2}=\infty $. The aim of this paper is to study
this case, to point out several new recent techniques and ideas and to
comment on the order of magnitude of the probabilities of large deviations.
Our proofs will be based on estimations of the Laplace transform.

One of our results is that for a strongly mixing sequence of centered and
bounded random variables satisfying, for a certain $c>0$, 
\begin{equation}
\alpha (n)\leq \exp (-2cn) \, ,  \label{alphacond}
\end{equation}%
we can find two constants $c_{1}$ and $c_{2}$ depending only on $c$ and on the uniform
bound of the random variables, such that, for all $x>0$, 
\begin{equation*}
\mathbb{P}(|S_{n}|>x)\leq \exp (-c_{1}x^{2}/n)+\exp (-c_{2}x/(\log n)(\log
\log n)) \, .
\end{equation*}%
Then, we use this exponential inequality and the techniques that lead to
this result to obtain moderate deviations asymptotic results, improving
Proposition 2.4 in Merlev\`{e}de and Peligrad (2009). Our results show that
we can come close up to a logarithmic term to the moderate deviations
asymptotics for independent random variables. Of course a kind of correction
is needed since the traditional large deviations results do not hold for
geometrically strongly mixing sequences. As a matter of fact the large
deviations principle does not hold even in the context of uniformly mixing
sequences with exponential rates. See Bryc and Dembo (1996), Example 1,
Proposition 5 and Example 2, that point out examples of empirical processes
of Doeblin recurrent Markov chains, that are therefore $\phi$-mixing with
exponential mixing rate and that do not satisfy the large deviations
principle.

We mention that strongly mixing sequences form a larger class than absolutely
regular sequences. The problem of the moderate deviations principle was
studied for absolutely regular Markov chains with exponential rates in de
Acosta (1997) and also in Chen and de Acosta (1998), when the transition
probabilities are stationary and there is a certain restriction on the class
of initial distributions. A class of processes satisfying a splitting
condition closely related to absolutely regular processes was considered by
Tsirelson (2008). Recently, Dedecker, Merlev\`{e}de, Peligrad and Utev
(2009) considered projective conditions with applications to $\phi$- mixing.

Notice that we do not require any degree of stationarity for obtaining
Bernstein inequality except for a uniform bound for the variables.

The strong mixing coefficient used in this paper can be generalized by using
smaller classes of functions than those used in Definition (\ref{defalpha2})
to include even more examples. Such examples include function of linear
processes with absolutely regular innovations and Arch models. In this paper,
we give an application to the moderate deviations principle for Kernel
estimators of the common marginal density of a certain class of continuous
time processes.

For the clarity of the proofs it will be more convenient to embed the
initial sequence into a continuous time process; namely, $(X_{t},t\geq 0)$
is defined from the original sequence $(X_{n},n\geq 1)$ by $X_{t}=X_{[t+1]}$%
. For a Borel set $A$, define 
\begin{equation}
S_{A}=\int\nolimits_{A}X_{t}dt \, .  \label{not}
\end{equation}%
Then $S_{[0,n]}=\sum\nolimits_{k=1}^{n}X_{k}$. The strong mixing coefficient
of $(X_{t},t\geq 0)$ is defined as 
\begin{equation*}
\tilde \alpha (u):=\sup_{t\geq 0}\alpha \left( \mathcal{M}_{t},\mathcal{G}%
_{t+u}\right) \rightarrow 0\mbox{ as }n\rightarrow \infty \, ,
\end{equation*}%
where $\mathcal{M}_{t}:=\sigma (X_{v},v\leq t)$ and $\mathcal{G}_{w}:=\sigma
(X_{v},v\geq w)$ .

\medskip

Notice that, since $\tilde \alpha (u) \leq \alpha ([u])$, if $(X_n, n \geq 1)
$ satisfies (\ref{alphacond}), the continuous type mixing coefficients still
satisfy a geometrically mixing condition; namely for any $u \geq 2$, 
\begin{equation}
\tilde \alpha (u) \leq \exp (-cu) \, .  \label{alphacondc}
\end{equation}

In the rest of the paper, $\|Y\|_{\infty}$ stands for the essential supremum of a random variable $Y$.

\section{Results}

Our first result is the following exponential inequality:

\begin{theorem}
\label{Bernstein}Let $(X_{j})_{j\geq 1}$ be a sequence of centered real-valued 
random variables. Suppose that the sequence satisfies (\ref{alphacond}%
) and that there exists a positive $M$ such that $\sup_{i\geq 1}\Vert
X_{i}\Vert _{\infty }\leq M$. Then there are positive constants $C_1$ and $C_2
$ depending only on $c$ such that for all $n\geq 4$ and $t$ satisfying $%
\displaystyle
0<t<\frac{1}{C_1 M(\log n)(\log \log n)}$, we have%
\begin{equation*}
\log \mathbb{E}\big (\exp (tS_{n}) \big )\leq \frac{C_2t^{2}nM^{2}}{1-C_1tM(\log n)(\log
\log n)}\text{ ,}
\end{equation*}%
In terms of probabilities, there is a constant $C_3$ depending on $c$ such
that for all $n\geq 4$ and $x\geq 0$ 
\begin{equation}
\mathbb{P}(|S_{n}|\geq x)\leq \exp \big (-{\frac{C_3x^{2}}{nM^{2}+Mx(\log n)(\log
\log n)}} \big )\, .  \label{bern1}
\end{equation}
\end{theorem}

As a counterpart, the following Bernstein type inequality holds.

\begin{theorem}
\label{bernstein2}Under conditions of Theorem \ref{Bernstein}, there are
positive constants $C_1$ and $C_2$ depending only on $c$ such that for all $%
n\geq 2$ and any positive $t$ such that $t<\frac{1}{C_1 M(\log n)^2}$, the
following inequality holds: 
\begin{equation*}
\log \mathbb{E}\big (\exp (tS_{n}) \big )\leq \frac{C_2t^{2}(nv^{2}+M^{2})}{%
1-C_1tM(\log n)^{2}}\,,
\end{equation*}%
where $v^{2}$ is defined by 
\begin{equation}
v^{2}=\sup_{i>0}\Bigl (\mathrm{Var}(X_{i})+2\sum_{j>i}|\mathrm{Cov}%
(X_{i},X_{j})|\Bigr)\,.  \label{condvn1}
\end{equation}%
In terms of probabilities, there is a constant $C_3$ depending only on $c$
such that for all $n\geq 2$, 
\begin{equation}
\mathbb{P}(|S_{n}|\geq x)\leq \exp \big (-{\frac{C_3x^{2}}{v^{2}n+M^{2}+xM(\log
n)^{2}}} \big )\,.  \label{bern2}
\end{equation}
\end{theorem}

To compare these two results, we notice that the coefficient of $x$ in the
inequality (\ref{bern1}) has a smaller order than the corresponding in (\ref%
{bern2}). However, the term $v^{2}n$ can be considerably smaller than $nM^{2}
$, which is an advantage in some applications when the variables are not
uniformly bounded. Notice also that if stationarity is assumed, $v^{2} $ can
be taken as 
\begin{equation*}
v^{2}=\mathrm{Var}(X_{1})+4\sum\nolimits_{i\geq 1}{\mathbb{E}}%
X_{1}^{2}I(|X_{1}|\geq Q(2\alpha _{i}))\, ,
\end{equation*}
where $Q(u)=\inf \{t>0,\mathbb{P}(|X_1|>t)\leq u\}$ for $u$ in $]0,1]$.

In the context of bounded functions $f$  of stationary geometrically strongly
mixing Markov chains, Theorem 6 in Adamczak (2008) provides a Bernstein type
inequality for $S_n(f) = f(X_1) +\dots + f(X_n)$ with the factor $\log n$ instead of $(\log n)^2$, which appears in
(\ref{bern2}). To be more precise, under the centering condition  $\mathbb{E} (f(X_1))=0$, he proves that
$$
\mathbb{P}(|S_{n}(f)|\geq x)\leq C\exp \Big (-\frac{1}{C} \min \Big (\frac{x^2}{n \sigma^2}, \frac{x}{\log n} \Big) \Big ) \, ,
$$
where $\sigma^2 = \lim_n n^{-1} {\rm Var} S_n(f)$.

\bigskip

The two previous results are useful to study the moderate deviations
principle (MDP) for the partial sums of the underlying sequences. In our
terminology the moderate deviations principle (MDP) signifies the following
type of behavior.

\begin{definition}
We say that the $MDP$ holds for a sequence $(T_n)_n$ of random variables
with the speed $a_{n}\rightarrow 0$ and rate function $I(t)$ if for each Borel set $A$,%
\begin{eqnarray}
-\inf_{t\in A^{o}}I(t) &\leq &\lim \inf_{n}a_{n}\log {\mathbb{P}}( \sqrt{%
a_{n}}T_{n}\in A)  \notag \\
&\leq &\lim \sup_{n}a_{n}\log {\mathbb{P}}(\sqrt{a_{n}}T_{n}\in A)\leq
-\inf_{t\in \bar{A}}I(t)\,  \label{mdpdef}
\end{eqnarray}%
where $\bar{A}$ denotes the closure of $A$ and $A^{o}$ the interior of $A$.
\end{definition}

Notice that the moderate deviations principle for $\big
(S_n/\sqrt{n} \big )$ is an intermediate behavior between CLT,\ ${\mathbb{P}}%
(S_{n}/\sqrt{n}\in A)$ and large deviation, ${\mathbf{P}}(S_{n}/n\in A)$.
Our moderate deviations results are the following:

\begin{theorem}
\label{thmMDP} Let $(X_{j})_{j\geq 1}$ be a sequence of centered real valued
random variables satisfying the assumptions of Theorem \ref{Bernstein}. Let $%
S_{n}=\sum_{i=1}^{n}X_{i}$, $\sigma _{n}^{2}=\mathrm{Var}(S_{n})$ and assume
in addition that $\lim \inf_{n\rightarrow \infty }\sigma _{n}^{2}/n>0$. Then
for all positive sequences $a_{n}$ with%
\begin{equation}
a_{n}\rightarrow 0\text{ and }\frac{na_{n}}{(\log n)^{2}(\log \log n)^{2}}%
\rightarrow \infty  \label{condan1}
\end{equation}%
the sequence $(\sigma _{n}^{-1}S_{n})_{n\geq 1}$ satisfies (\ref{mdpdef})
with the good rate function $I(t)=t^{2}/2$.
\end{theorem}

If we assume that the sequence is ${\mathbb{L}}_{2}$-stationary, then by
Lemma 1 in Bradley (1997), we get the following corollary:

\begin{corollary}
Let $(X_{j})_{j\geq 1}$ be as in Theorem \ref{thmMDP}. Suppose in addition
that the sequence is $\mathbb{L}_{2}-$stationary and $\sigma
_{n}^{2}\rightarrow \infty $. Then, $\lim_{n\rightarrow \infty }\sigma
_{n}^{2}/n=\sigma ^{2}>0$ and for all positive sequences $a_{n}$ satisfying (%
\ref{condan1}), $(n^{-1/2}S_{n})_{n\geq 1}$ satisfies (\ref{mdpdef}) with
the good rate function $I(t)=t^{2}/(2\sigma ^{2})$.
\end{corollary}

In the next result, we derive conditions ensuring that the MDP holds for the
partial sums of triangular arrays of strongly mixing sequences. For a double
indexed sequence $(X_{j,n},j\geq 1)_{n\geq 1}$ of real valued random
variables, we define for any $k\geq 0$, 
\begin{equation}
\alpha _{n}(k)=\sup_{j\geq 1}\alpha (\sigma (X_{n,i},i\leq j),\sigma
(X_{n,i},i\geq k+j)\,.  \label{defalphatrian}
\end{equation}

\begin{theorem}
\label{thmMDPtrian} For all $n\geq 1$, let $(X_{j,n},j\geq 1)_{n\geq 1}$ be
a double indexed sequence of centered real valued random variables such that
for every $j\geq 1$ and every $n\geq 1$, $\Vert X_{j,n}\Vert _{\infty }\leq
M_{n}$ where $M_{n}$ is a positive number. For all $n\geq 1$ and all $k\geq 0
$, let $\alpha _{n}(k)$ be defined by (\ref{defalphatrian}) and assume that $%
\alpha (k)=\sup_{n\geq 1}\alpha _{n}(k)$ satisfies (\ref{alphacond}). Define 
$v^{2}$ by 
\begin{equation}
v^{2}=\sup_{n\geq 1}\sup_{i>0}\Bigl (\mathrm{Var}(X_{i,n})+2\sum_{j>i}|%
\mathrm{Cov}(X_{i,n},X_{j,n})|\Bigr)\,.  \label{condvn}
\end{equation}%
and suppose $v^{2}<\infty .$ Let $S_{n}=\sum_{i=1}^{n}X_{i,n}$ , $\sigma
_{n}^{2}=\mathrm{Var}(S_{n})$ and assume in addition that $\lim
\inf_{n\rightarrow \infty }\sigma _{n}^{2}/n>0$. Then for all positive
sequences $a_{n}$ with%
\begin{equation}
a_{n}\rightarrow 0\text{ and }na_{n}/M_{n}^{2}(\log n)^{4}\rightarrow \infty
\label{condan2}
\end{equation}%
the sequence $(\sigma _{n}^{-1}S_{n})_{n\geq 1}$ satisfies (\ref{mdpdef})
with the good rate function $I(t)=t^{2}/2$.
\end{theorem}

\section{Discussion and Examples}

\medskip

\textbf{1.} The first comment is on Theorems \ref{Bernstein} and \ref%
{bernstein2}. Notice that compared to the traditional Bernstein inequality
for independent random variables there is a logarithmic correction in the
linear term in $x$ appearing in the inequalities (\ref{bern1}) and (\ref{bern2}). We includ in 
our paper another bound. Corollary \ref%
{bernstein3} gives better results than the other exponential bound results
in the large deviation range, that is when $x$ is close to $n$. As a matter
of fact the tail probability $\mathbb{P}(|S_{n}|\geq x)$ can be bounded with
the minimum of the right hand sides of inequalities (\ref{bern1}),(\ref%
{bern2}) and (\ref{bern3}). Among these inequalities, (\ref{bern1}) provides
the best condition leading to a moderate deviations principle when the
random variables are uniformly bounded.

\medskip

\textbf{2.} The strong mixing coefficients are not used in all their
strength. For obtaining our Bernstein type inequalities we can considerably
restrict the class of functions used to define the strong mixing
coefficients to those functions that are coordinatewise nondecreasing, and one
sided relations. Assume that for any index sets $Q$ and $Q^{\ast }$ (sets of
natural numbers) such that $Q\subset (0,p]$ and $Q^{\ast }\subset \lbrack
n+p,\infty )$, where $n$ and $p$ are integers, there exists a decreasing
sequence $\alpha ^{\ast }(n)$ such that 
\begin{equation*}
\mathrm{Cov}(f(S_{Q}),g(S_{Q^{\ast }}))\leq {\alpha }^{\ast
}(n)||f(S_{Q})||_{\infty }||f(S_{Q^{\ast }})||_{\infty }\,,
\end{equation*}%
where $f$ and $g$ are bounded functions coordinatewise nondecreasing. Here $%
S_{Q}=\sum\nolimits_{i\in Q}X_{i}$. Clearly the families of functions $\exp
(\sum\nolimits_{i\in Q}tx_{i})$ are coordinatewise nondecreasing for $t>0$
and then, for bounded random variables we have for all $t>0$ 
\begin{equation}
\mathrm{Cov}(\exp (tS_{Q}),\exp (tS_{Q^{\ast }}))\leq {\alpha }^{\ast
}(n)||\exp (tS_{Q})||_{\infty }||\exp (tS_{Q^{\ast }}))||_{\infty }\,.
\label{1}
\end{equation}%
Also by using the functions $f(x)=g(x)=x$%
\begin{equation}
\mathrm{Cov}(X_{j,}X_{j+n})\leq {\alpha }^{\ast }(n)||X_{j}||_{\infty
}||X_{j+n}||_{\infty }\,.  \label{2}
\end{equation}%
As a matter of fact these are the only functions we use in the proof of our
Bernstein inequality. So if ${\alpha }^{\ast }(n)$ decreases geometrically
our results still hold. Inequality (\ref{1}) is used to bound the Laplace
transform of partial sums, and Inequality (\ref{2}) is used to bound their
variance. Since both of these inequalities (\ref{1}) and (\ref{2}) are
stable under convolution we can obtain Bernstein type inequality for example
for sequences of the type $X_{n}=Y_{n}+Z_{n},$ where $Y_{n}$ is strongly
mixing as in Theorem \ref{Bernstein} and is $Z_{n}$ a noise, independent on $%
Y_{n}$, negatively associated, such as a truncated Gaussian sequence with
negative correlations.

We point out that similar results can be obtained by using alternative
mixing coefficients such as the $\tau $-mixing coefficient introduced by
Dedecker and Prieur (2004). Consequently we can treat all the examples in
Merlev\`{e}de, Peligrad and Rio (2009), namely: instantaneous functions of
absolutely regular processes, functions of linear processes with absolutely
regular innovations and ARCH($\infty )$ models.

\medskip

\textbf{3.} We now give an application to the moderate deviations principle
behavior for kernel estimators of the density of a continuous time process.

Let $X=(X_{t},t\geq 0)$ be a real valued continuous time process with an
unknown common marginal density $f$. We wish to estimate $f$ from the data $%
(X_{t},0\leq t\leq T)$. In what follows, we will call a kernel a
function $K$ from ${\mathbb{R}}$ to ${\mathbb{R}}$ which is a bounded
continuous symmetric density with respect to Lebesgue measure and such that 
\begin{equation*}
\lim_{|u|\rightarrow \infty }uK(u)=0,\text{ and }\,\int_{{\mathbb{R}}%
}u^{2}K(u)du<\infty \,.
\end{equation*}%
The kernel density estimator is defined as 
\begin{equation*}
f_{T}(x)=\frac{1}{Th_{T}}\int_{0}^{T}K\Big (\frac{x-X_{t}}{h_{T}}\Big )dt\,,
\end{equation*}%
where $h_{T}\rightarrow 0^{+}$ and $K$ is a kernel. 
In order to derive sufficient conditions ensuring that the MDP holds for the
sequence $\sqrt{T}(f_{T}(x)-{\mathbb{E}}f_{T}(x))$, we assume that there exists a constant $%
c>0$ such that for any $u\geq 0$, 
\begin{equation}
\alpha _{u}=\sup_{t\geq 0}\alpha (\sigma (X_{s},s\leq t),\sigma (X_{s},i\geq
u+t)\leq e^{-cu}\,.  \label{condalphacont*}
\end{equation}%
In addition, we assume that the joint distribution $f_{X_{s},X_{t}}$ between 
$X_{s}$ and $X_{t}$ exists and that $f_{X_{s},X_{t}}=f_{X_{0},X_{|t-s|}}$.
Applying Theorem \ref{thmMDPtrian}, we obtain the following result:

\begin{corollary}
\label{thmkernel} Suppose that $g_{u}=f_{X_{0},X_{u}}-f\otimes f$ exists for 
$u\neq 0$, and that the function $u\mapsto \sup_{x,y}\vert g_{u}(x,y)\vert$ 
is integrable on $]0,\infty \lbrack $ and $g_{u}(.,.)$ is continuous at $(x,x)$ 
for each $u>0$. In addition assume that the strong mixing coefficients of the process satisfy (%
\ref{condalphacont*}). Then for all positive sequences $a_{T}$ with%
\begin{equation*}
a_{T}\rightarrow 0,\text{\ }\frac{a_{[T]}}{a_{T}}\rightarrow 1\,,\text{ and }%
\frac{a_{T}Th_{T}^{2}}{(\log T)^{4}}\rightarrow \infty \text{ },
\end{equation*}%
the sequence $\sqrt{T}\big (f_{T}(x)-{\mathbb{E}}f_{T}(x) \big ) $ satisfies
(\ref{mdpdef}) with speed $a_{T}$ and the good rate function 
\begin{equation}
I(t)=t^{2}/\Big (4\int_{0}^{\infty }g_{u}(x,x)du\Big )^{-1}\,.
\label{rateker}
\end{equation}%
Furthermore if $f$ is differentiable and such that $f^{\prime }$ is $l$%
-Lipschitz for a positive constant $l$, and if $a_{T}Th_{T}^{4}\rightarrow 0$%
, then the sequence $\sqrt{T}(f_{T}(x)-f(x))$ satisfies (\ref{mdpdef}) with
speed $a_{T}$ and the good rate function defined by (\ref{rateker}).
\end{corollary}

Some examples of diffusion processes satisfying Condition (\ref%
{condalphacont*}) may be found in Veretennikov (1990) (see also Leblanc,
1997).

\section{Proofs}

First let us comment on the variance of partial sums. By using the notation (%
\ref{not}), for any compact set $K_A$ included in $[a, a +A]$ where $A>0$
and $a \geq 0$, we have that 
\begin{equation*}
\mathrm{Var} (S_{K_A}) \leq A \sup_{i>0}\Bigl (\mathrm{Var}%
(X_{i})+2\sum_{j>i}|\mathrm{Cov}(X_{i},X_{j})|\Bigr) \, . 
\end{equation*}
If the variables are bounded by $M$, then by using the definition (\ref%
{defalpha2}), we get that 
\begin{equation*}
\mathrm{Var} (S_{K_A}) \leq A\bigl(1+8\sum\nolimits_{i\geq 1}\alpha_{i}\bigr)M^{2}\leq K A M^{2} \, .
\end{equation*}%
If some degrees of stationarity are available we can have better upper
bounds. For instance if ${\mathbb{P}}(|X_{n}|>x)\leq {\mathbb{P}}(|X_{0}|>x),
$ then by Theorem 1.1 in Rio (2000), 
\begin{equation*}
\mathrm{Var} (S_{K_A}) \leq A \big (\mathrm{Var}(X_{0})+4\sum\nolimits_{i%
\geq 1}{\mathbb{E}}X_{0}^{2}I(|X_{0}|\geq Q(2\alpha _{i})) \big ) \, ,
\end{equation*}
where $Q(u)=\inf \{t>0,\mathbb{P}(|X_0|>t)\leq u\}$ for $u$ in $]0,1]$.

\subsection{Preliminary lemmas}

The first step is to prove an upper bound on the Laplace transform, valid
for small values of $t$. Without restricting the generality it is more
convenient to embed the index set into continuous time. In the following we
shall use the notation (\ref{not}).

\begin{lemma}
\label{propinter1} Let $(X_{n})_{n\geq 1}$ be as in Theorem \ref{Bernstein}.
Let $B \geq 2$ and $a \geq 0$. Then for any subset  $K_{B}$ of $(a,a+B]$
which is a finite union of intervals, and for any positive $t$ with $tM\leq %
\big ( \frac{1}{2} \big ) \wedge \big ( \frac{c}{2B} \big )^{1/2}$, we have 
\begin{equation}
\log \mathbb{E}\exp (tS_{K_{B}})\leq B \Bigl(6.2t^{2}v^{2}+ \frac{Mt}{2}\exp
\bigl(-\frac{c}{2tM}\bigr) \Bigr) \, ,  \label{resultpropinter}
\end{equation}%
where $v^{2}$ is defined by (\ref{condvn1}).
\end{lemma}

\begin{remark}
Notice that under our conditions $v^{2}\leq K M^{2}$ where $%
K=1+8\sum\nolimits_{i\geq 1}\alpha _{i}$.
\end{remark}

\noindent \textbf{Proof of Lemma \ref{propinter1}}. If $tM\leq 4/B$, then $%
tS_{K_{B}} \leq 4$, which ensures that 
\begin{equation*}
\exp (tS_{K_{B}})\leq 1+tS_{K_{B}}+{\frac{e^{4}-5}{16}}t^{2}(S_{K_{B}})^{2}%
\,,
\end{equation*}%
since the function $x\mapsto x^{-2}(e^{x}-x-1)$ is increasing. Now $%
(e^{4}-5)/16\leq 3.1$. Hence 
\begin{equation}
\mathbb{E}\exp (tS_{K_{B}})\leq 1+3.1Bv^{2}t^{2},  \label{pr}
\end{equation}
which implies Lemma \ref{propinter1} by taking into account that $\log
(1+x)\leq x$.

If $tM>4/B$, it will be convenient to apply Lemma \ref{lmaibra} in Appendix,
to get the result. Let $p$ be a positive real to be chosen later on. Let $%
k=[B/2p]$, square brackets denoting the integer part. We divide the interval 
$(a,a+B]$ into $2k$ consecutive intervals of equal size $B/(2k)$.
Denote these subintervals by $\{I_{j};1\leq j\leq 2k\}$ and let 
\begin{equation*}
\tilde S_{1}=\sum\limits_{j=1}^{k}S_{K_B \cap I_{2j-1}}\text{ and } \tilde
S_{2}=\sum\limits_{j=1}^{k}S_{K_B \cap I_{2j}}.
\end{equation*}%
By the Cauchy-Schwarz inequality, 
\begin{equation}
2\log \mathbb{E}\exp (tS_{K_{B}})\leq \log \mathbb{E}(2t \tilde S_{1})+\log 
\mathbb{E}\exp (2t \tilde S_{2}).  \label{pr2}
\end{equation}
Now let $p=1/(tM)$. Since $(X_{n})_{n\geq 1}$ satisfies Condition (\ref%
{alphacond}), and since $B/(2k) \geq p \geq 2$, by applying Lemma \ref%
{lmaibra} in Appendix, we obtain 
\begin{equation*}
\mathbb{E}\exp (2t \tilde S_{2})\leq k\exp \Bigl( \frac{MBt}{2}-\frac{cB}{2k}
\Bigr)+\prod_{j=1}^k\mathbb{E}\exp (2tS_{K_B \cap I_{2j}}).
\end{equation*}%
Notice that we are in the case $tM>4/B$ implying that $p\leq B/4$ and then $%
k\geq 2$. Now, under the assumptions of Lemma \ref{propinter1}, we have $%
tM\leq (c/(2B))^{1/2}$ which ensures that 
\begin{equation*}
MBt-\frac{cB}{2k}\leq MBt-cp\leq MBt-\frac{c}{tM}\leq -\frac{c}{2Mt}.
\end{equation*}%
Therefore, 
\begin{equation*}
\mathbb{E}(\exp (2t \tilde S_{2})\leq {\frac{BMt}{2}}\exp
(-c/(2tM))+\prod_{j=1}^{k}\mathbb{E}\exp (2tS_{K_B \cap I_{2j}}).
\end{equation*}
Since the random variables $(X_{i})_{i\geq 1}$ are centered, the Laplace
transforms of $\tilde S_{2}$ and each of $(S_{K_B \cap I_{2j}})_{j\geq 1}$ are
greater than one. Hence applying the inequality 
\begin{equation}
|\log x-\log y|\leq |x-y|\text{ for $x\geq 1$ and $y\geq 1$} \, ,
\label{eqlog}
\end{equation}%
we derive that 
\begin{equation*}
\log \mathbb{E}(\exp (2t \tilde S_{2}))\leq \sum_{j=1}^{k}\log \mathbb{E}%
\exp (2tS_{K_B \cap I_{2j}})+{\frac{BMt}{2}}\exp (-c/(2tM)).
\end{equation*}

Next $|2tS_{K_B \cap I_{2j}}|\leq 2tMB/(2k)$. Since $p\leq B/4$, $k \geq
B/(4p)$ implying that $|2tS_{K_B \cap I_{2j}}|\leq 4$, and consequently we
may repeat the arguments of the proof of (\ref{pr}), so that 
\begin{equation*}
\sum_{j=1}^{k}\log \mathbb{E}\exp (2tS_{K_B \cap I_{2j}})\leq 6.2Bt^{2}v^{2}.
\end{equation*}%
It follows that 
\begin{equation*}
\log \mathbb{E}\exp (2t \tilde S_{2})\leq 6.2Bt^{2}v^{2}+(BMt/2)\exp
(-c/(2tM)).
\end{equation*}%
Clearly the same inequality holds true for the log-Laplace transform of $%
\tilde S_{1}$ which, together with relation (\ref{pr2}), gives the result. $%
\diamond $

\bigskip

The key lemma for proving our theorems is a new factorization lemma. Its
proof combines the ideas of Bernstein big and small type argument with a
twist, diadic recurrence and Cantor set construction.

\begin{lemma}
\label{factor2}Let $(X_{i})_{i\geq 1}$ be as in Theorem \ref{Bernstein}.
Then, for every $A\geq 2 ( c \vee 10)$ there exists a subset $K_{A}$ of $%
[0,A]$, with Lebesgue measure larger than $A/2$ (not depending on the random
process)$\ $such that for all $t$, $0\leq tM\leq c_0/(\log A)$ where $c_0 = 
\frac{c}{8} \wedge \sqrt{\frac{c \log 2}{8}}$ 
\begin{equation}
\log (\mathbb{E}\exp (tS_{K_{A}}))\leq 6.2v^{2}t^{2}A+(c+1)A^{-1}\exp
(-c/4tM)\,.  \label{ineqexp2}
\end{equation}%
where $v^{2}$ is defined by (\ref{condvn1}). Moreover, if $A \geq 4 \vee (2c)
$ for all $0\leq tM< \frac{c \wedge 1}{2}$, we can find a constant $C$
depending only $c$ such that 
\begin{equation}
\log (\mathbb{E}\exp (tS_{(0,A]})\leq Ct^{2}M^2 A \log A \, .
\label{ineqexpA}
\end{equation}
\end{lemma}

\noindent \textbf{Proof of Lemma \ref{factor2}.} The proof is inspired by
the construction of a "Cantor set" and has several steps.

\medskip

\noindent \textbf{Step 1}. \textbf{A "Cantor set" construction}. Let $A$ be
a strictly positive real number strictly more than one. Let $\delta \in (0,1)$ that will be
selected later, and let $k_{A}$ be the largest integer $k$ such that $%
((1-\delta )/2)^{k}\geq 1/A.$ We divide the interval $[0,A]$ in three parts
and delete the middle one of size $A\delta .$ The remaining ordered sets are
denoted $K_{1,1},$ $K_{1,2}\ $\ and each has the Lebesgue measure $%
A(1-\delta )/2.$ We repeat the procedure. Each of the remaining two
intervals $K_{1,1},$ $K_{1,2},$ are divided in three parts and the central
one of length $A\delta (1-\delta )/2$ is deleted. After $j$ steps ($j\leq
k_{A}$), we are left with a disjoint union of $2^{j}$ intervals denoted by $%
K_{j,i}$, $1\leq i\leq 2^{j}$, each of size $A((1-\delta )/2)^{j}$ and we
deleted a total length $\sum\limits_{i=0}^{j-1}A\delta (1-\delta
)^{i}=A(1-(1-\delta )^{j}).$ We use the first index of sets $K_{j,i}$ to
denote the step, and second one to denote its order. Set $k=k_{A}$ when no
confusion is allowed, and define 
\begin{equation}
K_{A}=\bigcup_{i=1}^{2^{k}}K_{k,i}\,.  \label{defKA}
\end{equation}%
We shall use also the following notation: for any $\ell $ in $%
\{0,1,...,k_{A}\}$, 
\begin{equation*}
K_{A,\ell ,j}=\bigcup_{i=(j-1)2^{k-\ell }+1}^{j2^{k-\ell }}K_{k,i}\,,
\end{equation*}%
implying that for any $\ell $ in $\{0,1,...,k\}$: $K_{A}=\bigcup_{j=1}^{2^{%
\ell }}K_{A,\ell ,j}$.

\medskip

\noindent \textbf{Step 2. Proof of Inequality (\ref{ineqexp2}). } Here we
consider $K_{A}$ as constructed in step 1, with 
\begin{equation*}
\delta =\frac{\log 2}{2\log A}\,.
\end{equation*}%
Since $A\geq 2$, with this selection of $\delta $ we get that $\delta \leq
1/2$. Since $k(A)\leq \log A/\log 2$, it follows that 
\begin{equation*}
\lambda \big (\lbrack 0,A]\setminus K_{A}\big)\leq A\delta k(A)\leq A/2\text{
whence }\lambda (K_{A})\geq A/2\,.
\end{equation*}%
We estimate now the Laplace transform of $S_{K_{A}}$. We first  notice that
since $K_{A}$ is included in $[0,A]$, then if $tM\leq \big ( \frac{c_0}{
\log A}\big ) \wedge \sqrt{ \frac{c}{2A}}$, by applying Lemma \ref%
{propinter1}, we derive that 
\begin{equation*}
\log {\mathbb{E}} \exp (tS_{K_A}) \leq 6.2 A t^2 v^2 + A \frac{tM}{2}
\exp(-c/(2tM)) \, . 
\end{equation*}
Since $tM\leq c/(8\log A)$, we have 
\begin{equation}  \label{b1exp}
\exp (-c/(2tM))\leq A^{-2}\exp (-c/(4tM)) \, .
\end{equation}
Consequently
\begin{equation*}
\log {\mathbb{E}} \exp (tS_{K_A}) \leq 6.2 A t^2 v^2 + A^{-1} \frac{tM}{2}
\exp(-c/(4tM)) \, , 
\end{equation*}
proving (\ref{ineqexp2}) since $tM \leq 1/2$. Then we assume in the rest of
the proof that $(c/(2A))^{1/2}<tM\leq c_0/(\log A) $, and we shall then
estimate the Laplace transform of $S_{K_{A}}$ by the diadic recurrence. Let $%
t$ be a positive real. Since $K_{A,1,1}$ and $K_{A,1,2}$ are spaced by an
interval of size $A\delta $ and $A\delta \geq 2$ (since $A \geq 20$), by
using Lemma \ref{lmaibra} below and condition (\ref{alphacond}), we derive that 
\begin{eqnarray*}
\mathbb{E}\exp (tS_{K_{A}}) &=&\mathbb{E}\exp (tS_{K_{A,1,1}})\exp
(tS_{K_{A,1,2}}) \\
&\leq &\mathbb{E}\exp (tS_{K_{A,1,1}})\mathbb{E}\exp (tS_{K_{A,1,2}})+\exp
(-cA\delta +A(1 - \delta)tM)\,.
\end{eqnarray*}%
Since the variables are centered, $\mathbb{E}\exp (tS_{K_{A,1,i}})\geq 1$
for $i=1,2$. Hence by taking into account (\ref{eqlog}), we obtain that 
\begin{equation}  \label{dec1}
\log \mathbb{E}\exp (tS_{K_{A}})\leq \sum\nolimits_{i=1}^{2}\log \mathbb{E}%
\exp (tS_{K_{A,1,i}})+\exp (-cA\delta +A(1 - \delta)tM) \, .
\end{equation}
Now, let%
\begin{equation}
\ell =\ell (t)=\inf \{k\in \mathbb{Z}:A((1-\delta )/2)^{k}\leq \frac{c}{%
2(tM)^{2}}\}.  \label{defk}
\end{equation}%
Notice that $\ell (t)\geq 1$ since $t^{2}M^{2}>c/(2A)$. In addition by the
selection of $k_{A}$ and since $\delta \leq 1/2$ and $tM\leq \sqrt{\frac{c}{8%
}}$, it follows that $\ell (t)\leq k_{A}$. Notice also, by the bound on 
$tM$ and since $A \geq 4$, we have 
\begin{equation*}
A\delta \frac{(1-\delta )^{\ell(t) - 1}}{2^{\ell(t) - 1}} > \frac{c \delta}{%
2(tM)^{2}} \geq 2 \, . 
\end{equation*}
Using the homogeneity properties of $K_{A}$, the decomposition (\ref{dec1})
and iterating until $\ell(t)$, we get that 
\begin{eqnarray}  \label{factprem}
\log \mathbb{E}\exp (tS_{K_{A}})& \leq & \sum_{j=1}^{2^{\ell }}\log \mathbb{E%
}\exp (tS_{K_{A,\ell ,j}})  \notag \\
&+ &\sum_{j=0}^{\ell-1}2^{j}\exp \Big ( -cA\delta \frac{(1-\delta )^{j}}{%
2^{j}}+ 2tM A \frac{(1 - \delta)^{j+1}}{2^{j+1}} \Big )\, .
\end{eqnarray}
Consequently, for any $t\leq c\delta /(2M)$, 
\begin{equation}
\log \mathbb{E}\exp (tS_{K_{A}})\leq \sum_{j=1}^{2^{\ell }}\log \mathbb{E}%
\exp (tS_{K_{A,\ell ,j}})+\sum_{j=0}^{\ell-1}2^{j}\exp \Big ( -\frac{%
cA\delta }{2}\frac{(1-\delta )^{j}}{2^{j}}\Big ) \, .  \label{fact}
\end{equation}%
Whence, since $2^{\ell (t)}\leq A$ and $tM\leq c\delta /2$ we obtain 
\begin{equation*}
\sum_{j=0}^{\ell(t) -1}2^{j}\exp (-\frac{cA\delta }{2}\frac{(1-\delta )^{j}}{%
2^{j}})\leq 2^{\ell(t) }\exp \big (-c^{2}\delta /(2tM)^{2})\leq A\exp
(-c/(2tM) \big ) \, .
\end{equation*}%
Now we estimate each of the terms $\mathbb{E}\exp (tS_{K_{A,\ell ,j}}).$ By
the definition of $\ell(t) $ the conditions of Lemma \ref{propinter1} are
satisfied for $S_{K_{A,\ell ,j}}$ with $B=A((1-\delta )/2)^{\ell(t) }$.
Consequently, 
\begin{equation*}
\log \mathbb{E}\exp (tS_{K_{A,\ell ,j}})\leq B\big (6.2v^{2}t^{2}+tM\exp
(-c/(2tM))\big )\,.
\end{equation*}%
Therefore, by using (\ref{b1exp}), we derive that 
\begin{eqnarray*}
\log \mathbb{E}\exp (tS_{K_A}) &\leq &6.2v^{2}t^{2}A+tMA^{-1}\exp
(-c/4tM)+A^{-1}\exp (-c/(4tM)) \\
&\leq &6.2v^{2}t^{2}A+(c+1)A^{-1}\exp (-c/4tM).
\end{eqnarray*}%
This ends the proof of Inequality (\ref{ineqexp2}).

\medskip

\noindent \textbf{Proof of Inequality (\ref{ineqexpA}).\ }The proof of this
part uses the same construction with the difference that we do not remove
the holes from the set and we use instead their upper bound. Once again if $%
tM\leq (c/(2A))^{1/2}$,  applying Lemma \ref{propinter1} together with the
fact that $\exp(-c/(2tM)) \leq 2tM/c$ and $A \geq 4$, we derive that 
\begin{equation*}
\log {\mathbb{E}} \exp (tS_{(0,A]}) \leq A \log A (6.2 t^2 v^2 + (tM)^2/c )
\, , 
\end{equation*}
Taking into account that $v^2 \leq KM^2$ with $K = 1 + 8 \sum_{i \geq
1}\alpha_i$, the inequality (\ref{ineqexpA}) holds true with $C \geq 6.2 K +
1/c$. Then we can assume without loss of generality in the rest of the proof
that $(c/(2A))^{1/2}<tM< \frac{c \wedge 1}{2} $. We start by selecting $%
\delta =2tM/c<1.$ For this $\delta $, we select $k_{A}$ as before and $\ell
=\ell (t)\ $ as in relation (\ref{defk}). At first stage we divide as before
the interval $[0,A]$ in $3$ parts, the central one having a Lebesgue measure 
$A\delta .$ Notice that $A\delta \geq 2$ since $tM > \sqrt{c/(2A)} \geq c/A$
by the fact that $A \geq 2c$. Consequently, since the variables are bounded
by $M$, by condition (\ref{alphacond}), 
\begin{eqnarray*}
\mathbb{E}\exp (tS_{(0,A]}) &\leq &[\mathbb{E}\exp (tS_{K_{A,1,1}})\exp
(tS_{K_{A,1,2}})] e^{tAM\delta} \\
&\leq &[\mathbb{E}\exp (tS_{K_{A,1,1}})\mathbb{E}\exp (tS_{K_{A,1,2}})+\exp
(-A\delta c+AtM)] e^{tAM\delta} \, .
\end{eqnarray*}%
Since the variables are centered, $\mathbb{E}\exp (tS_{K_{A,1,i}}))\geq 1$
for $i=1,2$. Hence applying (\ref{eqlog}) and recalling that $\delta =2tM/c$%
, we obtain 
\begin{equation*}
\log \mathbb{E}\exp (tS_{[0,A]})\leq \sum\nolimits_{i=1}^{2}\log \mathbb{E}%
\exp (tS_{K_{A,1,i}}))+\exp (-A\delta c/2)+tAM\delta \, .
\end{equation*}%
Then, we repeat the same procedure starting with $K_{A,1,1}$ and $K_{A,1,2}$%
, and after $\ell = \ell(t) $ iterations we obtain 
\begin{eqnarray*}
\log \mathbb{E}\exp (tS_{(0,A]}) &\leq &\sum_{i=1}^{2^{\ell }}\log \mathbb{E}%
\exp (tS_{K_{A,\ell ,i}}) \\
&&+\sum_{i=0}^{\ell -1}2^{i} \Big ( \exp \big (-\delta c\frac{A(1-\delta
)^{i}}{2^{i+1}}\big )+tM\delta \frac{A(1-\delta )^{i}}{2^{i}} \Big ) \, .
\end{eqnarray*}
The above computation is valid since by the definition of $\ell(t)$ 
\begin{equation*}
A\delta \frac{(1-\delta )^{\ell(t) - 1}}{2^{\ell(t) - 1}} > \frac{c \delta}{%
2(tM)^{2}} = \frac{2tM}{2(tM)^{2}} \geq 2 \, . 
\end{equation*}
By the above considerations and the selection of $\ell $, proceeding as in
the proof of Inequality (\ref{ineqexp2}), we obtain 
\begin{equation*}
\log \mathbb{E}\exp (tS_{(0,A]})\leq 6.2v^{2}t^{2}A+A \big (\frac{tM}{2} + 1 %
\big )\exp (-c/(2tM))+tM\delta A\sum_{i=0}^{\ell \ -1}(1-\delta )^{j}.
\end{equation*}%
Now notice that for the selection $\delta =2tM/c$ and since $\ell \leq k_A$,
we have 
\begin{equation*}
tAM\delta \ell \leq 2t^{2}M^{2}A \log A / ( c \log 2) \, .
\end{equation*}
Also since for any $x \geq 0$, $\exp ( x) \geq x \vee (x^2/2)$, we get 
\begin{equation*}
\exp (-c/(2tM))\leq (2tM/c) \wedge (8t^{2}M^{2}/c^{2}) \, .
\end{equation*}%
Overall%
\begin{equation*}
\log \mathbb{E}\exp (tS_{(0,A]})\leq 6.2v^{2}t^{2}A+t^{2}M^{2}A ( 1/c +
8/c^2) + 2t^{2}M^{2}A \log A / ( c \log 2) \, .
\end{equation*}%
By taking into account that $v^{2}\leq KM^{2}$ we obtain the desired result
with the constant 
\begin{equation}  \label{bC}
C=6.2K+( 1/c + 8/c^2)+ 2/ ( c \log 2) \, ,
\end{equation}
where $K=(1+8\sum\nolimits_{i\geq 1}\alpha _{i}) $. $\diamond$

\bigskip

To prepare for the proof of Theorem \ref{Bernstein} we shall reformulate the
conclusions of Lemma \ref{factor2} in an alternative form. Keeping the same notations
 as in Lemma \ref{factor2}, the following corollary holds.

\begin{corollary}
Let $(X_{i})_{i\geq 1}$ be as in Theorem \ref{Bernstein}. Assume 
that $A \geq 2( c \vee 10)$ and $0\leq tM\leq c_0/(\log A)$ with $c_0
= \frac{c}{8} \wedge \sqrt{\frac{c \log 2}{8}}$. Then, there is a constant $%
C^{\prime }$ depending only on $c$ such that 
\begin{equation}
\log \mathbb{E}(\exp (tS_{K_{A}}))\leq \frac{C^{\prime 2}A(v+M/A)^{2}}{%
1-t(\log A)/c_0} \, .  \label{inex}
\end{equation}%
Assume that $A \geq 2( c \vee 10)$ and $0\leq tM< (c\wedge 1) /2$, then for
the constant $C$ defined in (\ref{bC}), 
\begin{equation}
\log \mathbb{E}(\exp (tS_{(0,A]})\leq \frac{Ct^{2}AM^{2}\log A}{%
1-2tM/(c\wedge 1)} \, .  \label{inex2}
\end{equation}
\end{corollary}

Before proving Theorem \ref{Bernstein} we remark that the second part of the
above corollary already gives a bound on the tail probability with a
correction in the quadratic term in $x$.

\begin{corollary}
\label{bernstein3}Under conditions of Theorem \ref{Bernstein}, for all $%
n\geq 2(c \vee 2)$ and $x\geq 0,$%
\begin{equation}
\mathbb{P}(|S_{n}|\geq x)\leq \exp \big (-\frac{x^{2}}{n(\log
n)4CM^{2}+4Mx/(c\wedge 1)} \big ) \, ,  \label{bern3}
\end{equation}
where $C$ is defined in (\ref{bC}).
\end{corollary}

\subsection{Proof of Theorem \protect\ref{Bernstein}.}

\bigskip If $n \leq 16(c \vee 10)^2$ then for any positive $t$ such that $tM
< \frac{1}{4(c \vee 10)^2}$ we get that $|t S_n| \leq 4$. Hence as in the
beginning of the proof of Lemma \ref{propinter1}, we derive that 
\begin{equation*}
\log {\mathbb{E}} (\exp(tS_n)) \leq 3.1 \frac{ nv^2 t^2}{1 - 4 tM(c \vee 10)^2} \,
. 
\end{equation*}
We assume now that $n \geq 16(c \vee 10)^2$. Let us first introduce the
following notation: for any positive real $A$, let $K_A$ be the Cantor set
as defined in step 1 of the proof of Lemma \ref{propinter1}, let $\lambda
(K_A)$ be the Lebesgue measure of $K_A$, and let $F_A$ be the nondecreasing
and continuous function from $[0,A]$ onto $[0, A- \lambda (K_A)]$ defined by 
\begin{equation}  \label{P1prop3}
F_A (t) = \lambda ([0,t]\cap K_A^c) \ \hbox{ for any } t\in [0,A],
\end{equation}
where $K_A^c = [0,A]\setminus K_A$. Let $F_{A}^{-1}$ be the inverse function of $F_A$. Let $A_0 = n$. 
Define then the
real-valued process $(X^{(1)}_t)_t$ from $(X_t)_{t\in [0,A_0]}$ by 
\begin{equation*}
X^{(1)}_t = X_{F_{A_0}^{-1} (t)} \ \hbox{ for any } t\in [0,{A_0}-\lambda
(K_{A_0})]. 
\end{equation*}
Let $A_1 = {A_0}-\lambda (K_{A_0})$. Clearly, the random process $%
(X^{(1)}_t)_{t \in [0,A_1]}$ is uniformly bounded by $M$ and verifies (\ref%
{alphacondc}) with the same constant. We now define inductively the sequence 
$(A_i)_{i\geq 0}$ and the random processes $(X_t^{(i)})_{i\in [0,A_i]}$ as
follows. First $A_0 = n$ and $(X^{(0)}_t) = (X_t)$. And second for any
nonnegative integer $i$, $A_{i+1} = A_i - \lambda (K_{A_i})$ and, for any $t$
in $[0,A_{i+1}]$, 
\begin{equation}  \label{P2prop3}
X^{(i+1)}_t = X_{F_{A_i}^{-1} (t)}^{(i)} \, .
\end{equation}
Then, for any nonnegative integer $j$, the following decomposition holds 
\begin{equation}  \label{P3prop3}
\int_0^n X_u du = \sum_{i=0}^{j-1} \int_{K_{A_i}} X^{(i)}_u du +
\int_0^{A_j} X^{(j)}_u du.
\end{equation}
Let 
\begin{equation}  \label{P4prop3}
Y_i= \int_{K_{A_i}} X^{(i)}_u du \text{ for $0\leq i \leq j-1$ and }Z_j=
\int_0^{A_j} X^{(j)}_u du.
\end{equation}
Now set 
\begin{equation*}
L=L_n = \inf \{ j \in {\mathbb{N}}^* \, , \, A_j \leq n/(\log n) \} \, . 
\end{equation*}
Notice that, since $A_j \leq n/2^j$, 
\begin{equation}  \label{bL}
L \leq [(\log \log n)/(\log 2)]+1 \, .
\end{equation}
Also by the definition of $L$, $A_{L-1} \geq n /(\log n)$. Since $\log n
\leq 2\sqrt n$, it follows that 
\begin{equation*}
A_{L-1} \geq \sqrt{n} /2 \geq 2(c \vee 10) \, . 
\end{equation*}%
Hence, we can apply the inequality (\ref{inex}) to each $Y_j$ for all $0\leq
j\leq L_{n} -1$. Consequently for every $0\leq j\leq L_{n} -1$, and any
positive $t$ satisfying $tM< c_0/(\log (n/2^{j}))$, 
\begin{equation}
\log \mathbb{E}(\exp (tY_j))\leq \frac{C^{\prime
2}(v(n/2^{j})^{1/2}+(n/2^{j})^{-1/2}M)^{2}}{1-Mt(\log (n/2^{j}))/c_0} \, .
\label{inex3}
\end{equation}
To estimate $Z_L$, we first assume that $A_L \geq 2 (c \vee 2)$. Applying
Inequality (\ref{inex2}) we then obtain, for any positive $t$ such that $%
tM<(c \wedge 1) /2$, 
\begin{equation*}
\log \mathbb{E}(\exp (tZ_L))\leq \frac{Ct^{2}M^{2}n}{1-2tM/(c\wedge 1)} \, .
\end{equation*}%
To aggregate all the contributions, we now apply Lemma \ref{breta} of
Appendix with $\kappa_{i}=M(\log (n/2^{j}))/c_0$ and $\sigma
_{i}^{2}=C^{\prime j})^{1/2}+(n/2^{j})^{-1/2}M)^{2}$ for $0\leq i\leq L-1,$ $%
\sigma _{L}^{2}=CnM^{2}$ and $\kappa_{L}=2M/(c\wedge 1)$. Consequently, by (%
\ref{bL}), there exists $C_1$ depending only on $c$ such that 
\begin{equation*}
\sum\nolimits_{i=1}^{L}\kappa_{i}+\kappa_{L}\leq C_1 M (\log n)(\log \log n)
\, .
\end{equation*}%
Furthermore 
\begin{equation*}
\sum\nolimits_{i=0}^{L-1}\sigma _{i}+\sigma _{L}\leq \sqrt{C^{\prime }}( 4
vn^{1/2}+ 2 M((\log n)/n)^{1/2})+\sqrt{Cn}M \, .
\end{equation*}%
Hence by Lemma \ref{breta}, for any $n \geq 4$ and any positive $t <1/(MC_1
\log n(\log \log n))$ there exists $C_2$ depending only on $c$ such that 
\begin{eqnarray*}
\log \mathbb{E}(\exp (tS_{n})) \leq \frac{C_2nt^{2}M^{2}}{1-tM C_1(\log
n)(\log \log n)} \, ,
\end{eqnarray*}%
and the result follows. If $A_L \leq 2 (c \vee 2)$, it suffices to notice
that if $tM < 2/(c \vee 2)$, then $|t Z_L| \leq 4$. Hence as in the proof of
Lemma \ref{propinter1}, we derive that 
\begin{equation*}
\log {\mathbb{E}} (\exp(tZ_L) )\leq 3.1 \frac{ 4 (c \vee 2)^2 t^2M^2}{1 - tM(c
\vee 2)/2} \, , 
\end{equation*}
and we proceed as before with $\kappa_{L} =M(c \vee 2)/2$ and $\sigma_L =
2(c \vee 2) M$. Inequality (\ref{bern1}) follows from the Laplace transform
estimate by standard computations. $\circ $

\subsection{Proof of Theorem \protect\ref{bernstein2}.}

We proceed as in the proof of Theorem \ref{Bernstein} with the  difference
that for $n \geq 2 ( c \vee 10)$, we choose 
\begin{equation*}
L=L_n = \inf \{ j \in {\mathbb{N}}^* \, , \, A_j \leq 2(c \vee 10) \} \, . 
\end{equation*}
Consequently, 
\begin{equation*}
L \leq \Big [ \frac{\log (n) - \log (2(c \vee 10))}{\log 2} \Big ] +1 \, .
\end{equation*}

\subsection{ Proof of Theorem \protect\ref{thmMDP}.}

The proof is based on the construction of the Cantor-like sets as described
in the proof of Lemma \ref{factor2}. Let $(\varepsilon _{n})_{n\geq 1}$ be a
sequence converging to $0$ that will be constructed later. Without loss of
generality we assume $\varepsilon _{n}<\log 2$ and define 
\begin{equation}
\delta _{n}=\frac{\varepsilon _{n}}{\log n} \, .  \label{selectdel1}
\end{equation}
We impose for the moment that $\varepsilon_n$ has to satisfy 
\begin{equation}  \label{restepsi1}
\delta_n \sqrt{n a_n} \rightarrow \infty \, .
\end{equation}%
(It is always possible to choose  such an $\varepsilon_n$ since (\ref%
{condan1}) is assumed). Select in addition 
\begin{equation*}
k_{n}=\inf \big \{j\in \mathbb{N}^{\ast }:n\frac{(1-\delta _{n})^{j}}{2^{j}}%
\leq \sqrt{na_{n}}\big \}.
\end{equation*}%
Construct the intervals $K_{k_{n},i}$, $1\leq i\leq 2^{k_n}$ , as in the
step 1 of the proof of Lemma \ref{factor2}. Then $R_{n}=(0,n]\backslash K_{n}
$, with $K_{n}=\cup _{i=1}^{2^{k_{n}}}K_{k_{n},i}$, has a Lebesgue measure
smaller than $\delta _{n}nk_{n}=o(n)\ $ and by Inequality (\ref{bern1})
there exists a constant $C$ depending only on $c$ such that 
\begin{equation*}
a_{n}\log \mathbb{P}(|S_{R_{n}}|\geq x\sigma _{n}/\sqrt{a_{n}})\leq -{\frac{%
Cx^{2}\sigma _{n}^{2}}{\delta _{n}k_{n}nM^{2}+Mx\sqrt{\sigma _{n}^{2}/a_{n}}%
(\log n)(\log \log n)}}\,.
\end{equation*}%
Taking into account that $\lim \inf_{n\rightarrow \infty }\sigma
_{n}^{2}/n>0 $, Condition (\ref{condan1}) ensures that 
\begin{equation*}
\lim_{n\rightarrow \infty }a_{n}\log \mathbb{P}(|S_{R_{n}}|\geq x\sigma _{n}/%
\sqrt{a_{n}})=-\infty .
\end{equation*}%
According to Theorem 4.2.13 in Dembo and Zeitouni (1998), $S_{R_{n}}$ is
negligible for the moderate deviations type of behavior. To treat the main
part we rewrite the inequality (\ref{factprem}) by using both sides of Lemma %
\ref{lmaibra} from Appendix and the fact that for $n$ large enough, (\ref%
{restepsi1}) entails that for any $t$ 
\begin{equation*}
\frac{|tM|}{\sigma_n \sqrt{a_n}} \leq \delta_n/2 \, \text{ and } \, n
\delta_n \frac{(1 - \delta_n)^{k_n - 1 }}{2^{k_n - 1 }} > \delta_n \sqrt{n
a_n} \geq 2 \, , 
\end{equation*}%
(by using also the fact that $\liminf_n \sigma_n^2/n >0$). So, by 
definition of the sets $K_{k_{n},i}$ from relation (\ref{defKA}), we get for any
real $t$ and $n$ large enough, 
\begin{equation*}
\Big |\log \mathbb{E}\exp \Big (\frac{tS_{K_{n}}}{\sigma _{n}\sqrt{a_{n}}}%
\Big )-\sum_{i=1}^{2^{k_{n}}}\log \mathbb{E}\exp \Big (\frac{tS_{K_{k_{n},i}}}{\sigma
_{n}\sqrt{a_{n}}} \Big )\Big |\leq \sum_{j=0}^{k_{n}-1}2^{j}\exp (-\frac{cn\delta _{n}%
}{2}\frac{(1-\delta _{n})^{j}}{2^{j}}).
\end{equation*}%
Now, by the definition of $k_{n}$ notice that $2^{k_{n}-1}\leq \sqrt{n/a_{n}}
$. Whence 
\begin{gather*}
\sum_{j=0}^{k_{n}-1}2^{j}\exp (-n\delta _{n}\frac{(1-\delta _{n})^{j}}{2^{j}}%
)\leq 2^{k_{n}}\exp (-\frac{c}{2}\delta _{n}\sqrt{na_{n}})\leq 2\sqrt{\frac{n%
}{a_{n}}}\exp (-\frac{c}{2}\delta _{n}\sqrt{na_{n}}) \\
\leq 2\frac{\sqrt{na_{n}}}{a_{n}}\exp (-\delta _{n}\sqrt{na_{n}})\leq \frac{2%
}{a_{n}}\exp (\log (\sqrt{na_{n}})-\frac{c}{2}\delta _{n}\sqrt{na_{n}}).
\end{gather*}%
We select now the sequence $\varepsilon _{n}$ used in the construction of $%
\delta _{n}$, in such a way that 
\begin{equation}  \label{selecteps}
\log (\sqrt{na_{n}})-\frac{c}{2}\delta _{n}\sqrt{na_{n}} \rightarrow \infty
\, .
\end{equation}%
Denote 
\begin{equation*}
A_{n}=\frac{\sqrt{na_{n}}}{(\log n)^{2}}\text{ and }B_{n}=\frac{\sqrt{na_{n}}%
}{(\log n)(\log \log n)}.
\end{equation*}%
If $na_{n}\geq (\log n)^{5}$ we select $\varepsilon _{n}=A_{n}^{-1/2}.$ If $%
na_{n}<(\log n)^{5}$ we take $\varepsilon _{n}=B_{n}^{-1/2}$. Notice that by
construction and by (\ref{condan1}), $\varepsilon _{n}\rightarrow 0$ as $%
n\rightarrow \infty $ and (\ref{restepsi1}) is satisfied. It is easy to see
that 
\begin{equation*}
\log (\sqrt{na_{n}})-c\delta _{n}\sqrt{na_{n}}/2\rightarrow -\infty \text{
when }n\rightarrow \infty .
\end{equation*}
Indeed, if $na_{n}\geq (\log n)^{5}$, then $\varepsilon _{n}=A_{n}^{-1/2}$,
so that 
\begin{eqnarray*}
\log (\sqrt{na_{n}})-c\delta _{n}\sqrt{na_{n}}/2 &\leq &(\log n)/2-\frac{c}{%
\log n}\frac{A_{n}}{2A_{n}^{1/2}}(\log n)^{2} \\
&=&(\log n)(1-cA_{n}^{1/2})/2 \, .
\end{eqnarray*}
If $na_{n}<(\log n)^{5}$, then $\varepsilon _{n}=B_{n}^{-1/2}$, so that 
\begin{gather*}
\log (\sqrt{na_{n}})-c\delta _{n}\sqrt{na_{n}}/2=\log (\frac{\sqrt{na_{n}}}{%
(\log n)^{5/2}})+\log [(\log n)^{5/2}] \\
-c\varepsilon _{n}B_{n}(\log \log n)/2\leq (5/2)\log \log
n-cB_{n}^{1/2}(\log \log n)/2 \, .
\end{gather*}%
Consequently, for any real $t$, we get 
\begin{equation*}
a_{n}\Big |\log \mathbb{E}\exp \Big (\frac{tS_{K_{n}}}{\sigma _{n}\sqrt{a_{n}}}%
\Big )-\sum_{i=1}^{2^{k_{n}}}\log \mathbb{E}\exp \Big (\frac{tS_{K_{k_{n},i}}}{\sigma
_{n}\sqrt{a_{n}}} \Big )\Big |\rightarrow 0\text{ as }n\rightarrow \infty \, .
\end{equation*}

Therefore the proof is reduced to proving the MDP for a triangular array of
independent random variables $S^*_{K_{k_{n},i}}$, $1\leq i\leq 2^{k_{n}}$,
each having the same law as $S_{K_{k_{n},i}}$. By the selection of $k_{n}$,
for any $1\leq i\leq 2^{k_{n}}$, $\Vert S_{K_{k_{n},i}}\Vert _{\infty }\leq M%
\sqrt{na_{n}}$. In addition, for each $\epsilon >0$%
\begin{equation*}
\lim_{n\rightarrow \infty }\frac{1}{\sigma _{n}^{2}}\sum%
\nolimits_{j=1}^{2^{k_{n}}}E(S_{K_{k_{n},j}}^{2}I(|S_{K_{k_{n},j}}|>\epsilon
\sigma _{n}\sqrt{a_{n}}))=0\,,
\end{equation*}%
by using Inequality (\ref{bern1}) to give an upper bound of ${\mathbb{P}}%
(|S_{K_{k_{n},i}}|\geq x)$. Hence, by Lemma \ref{arcones}, we just have to
prove that 
\begin{equation}
\lim_{n\rightarrow \infty }\frac{1}{\sigma _{n}^{2}}\sum_{j=1}^{2^{k_{n}}}{%
\mathbb{E}}(S_{K_{k_{n},j}})^{2}=1\,.  \label{cond1arc}
\end{equation}%
With this aim we first notice that since $\mathrm{Card}R_{n}=o(n)$, it
follows that when $n\rightarrow
\infty $, $\mathrm{Var}(S_{R_{n}})/n\rightarrow 0$. Hence to prove (\ref{cond1arc}), it suffices to prove that 
\begin{equation}
\mathrm{Var}(\sum_{i=1}^{2^{k_{n}}}S_{K_{k_{n},i}})/\sum_{i=1}^{2^{k_{n}}}%
\mathrm{Var}(S_{K_{k_{n},i}})\rightarrow 1\text{ as }n\rightarrow \infty \,.
\label{equivar}
\end{equation}%
Between the two sets $K_{k_{n},i}$ and $K_{k_{n},j}$, for $i\neq j$, there is a
gap at least equal to $n\delta _{n}2^{1-k_{n}}(1-\delta _{n})^{k_{n}-1} >
\delta _{n}\sqrt{na_{n}}$, by the selection of $k_n$. Consequently, since
for all $j$, $\Vert S_{K_{k_{n},j}}\Vert _{\infty }\leq M\sqrt{na_{n}}$, by (%
\ref{defalpha2}) we get 
\begin{equation}
\sum\nolimits_{i=1}^{2^{k_{n}}-1}\sum\nolimits_{j=i+1}^{2^{k_{n}}-1}\mathrm{%
cov}(S_{K_{k_{n},i}},S_{K_{k_{n},j}})\leq 4\times
2^{k_{n}}M^{2}na_{n}\sum_{j\geq 1}\exp (-cj\delta _{n}\sqrt{na_{n}})\,.
\label{covrio}
\end{equation}%
Since $2^{k_{n}}\leq 2\sqrt{n/a_{n}}$, by using (\ref{selecteps}), it
follows that 
\begin{equation}
\sum\nolimits_{i=1}^{2^{k_{n}}-1}\sum\nolimits_{j=i+1}^{2^{k_{n}}-1}\mathrm{%
Cov}(S_{K_{k_{n},i}},S_{K_{k_{n},j}})=o(n)\,,  \label{cov}
\end{equation}%
which together with the fact that $C_{1}n\leq \mathrm{Var}%
(\sum_{i=1}^{2^{k_{n}}}S_{K_{k,i}})\leq C_{2}n$ implies (\ref{equivar}).

\subsection{ Proof of Theorem \protect\ref{thmMDPtrian}.}

The proof is similar to that of Theorem \ref{thmMDP} with the following
modifications. Inequality (\ref{bern2}) is used instead of Inequality (\ref%
{bern1}) (notice that (\ref{condan2}) implies $n/M_{n}^{2}\rightarrow \infty 
$), and the sequences $\varepsilon _{n}$ (defining $\delta _{n}$) and $k_{n}$
are selected as follows: 
\begin{equation}
\varepsilon _{n}\rightarrow 0\,\text{ and }\,\varepsilon _{n}\frac{\sqrt{%
na_{n}}}{M_{n}(\log n)^{2}}\rightarrow \infty \,,  \label{selectvarep2}
\end{equation}%
and 
\begin{equation*}
k_{n}=\inf \big \{j\in \mathbb{N}^{\ast }:n(M_{n}\vee 1)\frac{(1-\delta
_{n})^{j}}{2^{j}}\leq \sqrt{na_{n}}\big \}.
\end{equation*}%
Notice that with this selection, $2^{k_{n}}\leq 2(M_{n}\vee 1)\sqrt{n/a_{n}}$%
. By the selection of $\varepsilon _{n}$, all the steps of the previous
theorem can be done similarly. Also to prove that 
\begin{equation*}
\mathrm{Var}(\sum_{i=1}^{2^{k_{n}}}S_{K_{k_{n},i}})/\sum_{i=1}^{2^{k_{n}}}%
\mathrm{Var}(S_{K_{k_{n},i}})\rightarrow 1\text{ as }n\rightarrow \infty \,,
\end{equation*}%
we make use of Condition (\ref{condvn}) together with the fact that by the
selection of $k_{n}$, for all $j$, $\Vert S_{K_{k_{n},j}}\Vert _{\infty
}\leq \sqrt{na_{n}}$. The inequality (\ref{covrio}) becomes 
\begin{eqnarray*}
&&\sum\nolimits_{i=1}^{2^{k_{n}}-1}\sum\nolimits_{j=i+1}^{2^{k_{n}}-1}%
\mathrm{Cov}(S_{K_{k_{n},i}},S_{K_{k_{n},j}}) \\
&&\quad \quad \leq 8(M_{n}\vee 1)na_{n}\sqrt{n/a_{n}}\sum_{j\geq 1}\exp
(-cj\delta _{n}(M_{n}\vee 1)^{-1}\sqrt{na_{n}})\,,
\end{eqnarray*}%
which implies (\ref{cov}) by the selection of $\delta _{n}$ and the fact
that $M_n=o(\sqrt n )$.

\subsection{\textbf{Proof of Corollary \protect\ref{thmkernel}.}}

For each $n\geq 1$, let us construct the following sequence of triangular
arrays: for any $i\in {\mathbb{Z}}$, 
\begin{equation*}
X_{n,i}=\frac{1}{\sqrt{\delta }h_{T}}\Big \{\int_{(i-1)\delta }^{i\delta }K%
\Big (\frac{x-X_{t}}{h_{T}}\Big )dt-{\mathbb{E}}\int_{(i-1)\delta }^{i\delta
}K\Big (\frac{x-X_{t}}{h_{T}}\Big )dt\Big \} \, ,
\end{equation*}%
where $n\delta =T$, $n=[T]$ , $(T\geq 1$) and consequently $2>\delta \geq 1$%
. Notice that 
\begin{equation*}
\sum_{i=1}^{n}X_{n,i}=T\big (f_{T}(x)-{\mathbb{E}}f_{T}(x) \big ) \, .
\end{equation*}%
Now for any $k\geq 1$, the strong mixing coefficients, $\alpha _{n}(k)$, of
the processes $(X_{n,i})_{i\in {\mathbb{Z}}}$ are uniformly bounded by the
strong mixing coefficient $\alpha _{k-1}$ of the process $(X_{t},t\in {%
\mathbb{R}})$. Hence to apply Theorem \ref{thmMDPtrian}, it suffices to show
(\ref{condvn}) and to prove that 
\begin{equation*}
T^{-1}\mathrm{Var}\big (\sum_{i=1}^{n}X_{n,i}\big )\rightarrow
2\int_{0}^{\infty }g_{u}(x,x)du \text{ as $n\rightarrow \infty $}\,.
\end{equation*}%
The above convergence was proved by Castellana and Leadbetter (1986) under
 assumptions on $g_{u}$. To prove (\ref{condvn}), we first notice
that for all $j>i$, 
\begin{equation*}
\mathrm{Cov}(X_{i,n},X_{j,n})=\frac{1}{\delta h_T^2}\int_{{\mathbb{R}}^{2}} K%
\Big (\frac{x-y}{h_{T}}\Big )K\Big (\frac{x-z}{h_{T}}\Big )\int_{i\delta - \delta
}^{i\delta }\int_{j\delta - \delta }^{j\delta }g_{t-s}(y,z)dsdtdy dz \,.
\end{equation*}%
Consequently, since $K$ is a kernel, for all $j>i$, 
\begin{equation*}
\big |\mathrm{Cov}(X_{i,n},X_{j,n})\big |\leq \int_{(j-i-1)\delta
}^{(j-i+1)\delta } \sup_{x,y}\vert g_{u}(x,y)\vert du\,.
\end{equation*}%
Similarly 
\begin{equation*}
\mathrm{Var}(X_{i,n})\leq 2\int_{0}^{\delta }\sup_{x,y}\vert g_{u}(x,y)\vert du\,.
\end{equation*}%
Hence (\ref{condvn}) holds with 
\begin{equation*}
v^{2}\leq 2\int_{0}^{\delta }\sup_{x,y}\vert g_{u}(x,y)\vert du+4\int_{0}^{\infty
}\sup_{x,y}\vert g_{u}(x,y)\vert du\,.
\end{equation*}%
To finish the proof, it remains to notice that if $f$ is differentiable and
such that $f^{\prime }$ is $l$-Lipschitz for a positive constant $l$ then,
since $K$ is a kernel, 
\begin{equation*}
|{\mathbb{E}}f_{T}(x)-f(x)|=O(h_{T}^{2})\,,
\end{equation*}%
(see for instance relation 4.15 in Bosq (1998)).

\section{\protect\Large Appendix.}

One of our tools is the technical lemma below, which provides bounds for the
log-Laplace transform of any sum of real-valued random variables. It comes
from Lemma 3 in Merlev\`{e}de, Peligrad and Rio (2009).

\begin{lemma}
\label{breta} Let $Z_{0},Z_{1},\ldots $ be a sequence of real valued random
variables. Assume that there exists positive constants $\sigma _{0},\sigma
_{1},\ldots $ and $\kappa_{0},\kappa_{1},\ldots $ such that, for any $i\geq 0
$ and any $t$ in $[0,1/c_{i}[$, 
\begin{equation*}
\log \mathbb{E}\exp (tZ_{i})\leq (\sigma _{i}t)^{2}/(1-\kappa_{i}t)\,.
\end{equation*}%
Then, for any positive $n$ and any $t$ in $[0,1/(\kappa_{0}+\kappa_{1}+%
\cdots +\kappa_{n})[$, 
\begin{equation*}
\log \mathbb{E}\exp (t(Z_{0}+Z_{1}+\cdots +Z_{n}))\leq (\sigma
t)^{2}/(1-\kappa t),
\end{equation*}%
where $\sigma =\sigma _{0}+\sigma _{1}+\cdots +\sigma _{n}$ and $\kappa =
\kappa_{0}+ \kappa_{1}+\cdots + \kappa_{n}$.
\end{lemma}

The next lemma is due to Arcones (Lemma 2.3, 2003) and it permits us to derive
the MDP for triangular array of independent r.v.'s.

\begin{lemma}[Arcones (2003)]
\label{arcones} Let $\{X_{n,j};1\leq j\leq k_{n}\}$ be a triangular array of
independent r.v.'s with mean zero. Let $\{a_{n}\}_{n\geq 1}$ be a sequence
of real numbers converging to $0$. Suppose that:\newline
(i) The following limit exists and is finite: 
\begin{equation*}
\lim_{n\rightarrow \infty }\sum\nolimits_{j=1}^{k_{n}}{\mathbb{E}}%
(X_{n,j}^{2})=\sigma ^{2}\,,
\end{equation*}%
(ii) There exists a constant $C$ such that for each $1\leq j\leq k_{n}$,%
\begin{equation*}
|X_{n,j}|\leq C\sqrt{a_{n}}\,,
\end{equation*}%
$\ $\linebreak (iii) For each $\epsilon >0$%
\begin{equation*}
\lim_{n\rightarrow \infty }\sum\nolimits_{j=1}^{k_{n}}{\mathbb{E}}%
(X_{n,j}^{2}I(|X_{n,j}|>\epsilon \sqrt{a_{n}})=0\,.
\end{equation*}%
Then for all real $t$, $a_{n}\sum_{j=1}^{k_{n}}\log {\mathbb{E}}\exp
(tX_{n,j})\rightarrow t^{2}\sigma ^{2}/2$ and consequently the MDP holds for 
$(\sum_{j=1}^{k_{n}}X_{j,n})$ with speed $a_{n}$ and good rate function $%
I(t)=t^{2}/(2\sigma ^{2})$.
\end{lemma}

We first recall the following lemma, which is a well-known corollary of
Ibragimov's covariance inequality for nonnegative and bounded random
variables.

\begin{lemma}[Ibragimov (1962)]
\label{lmaibra} Let $Z_{1}$, ..., $Z_{p}$ be real-valued nonnegative random
variables each a.s. bounded, and let 
\begin{equation*}
\alpha =\sup_{k\in \lbrack 1,p]}\alpha (\sigma (Z_{i}:i\leq k),\sigma
(Z_{i}:i>k))
\end{equation*}%
Then 
\begin{equation*}
\mathbb{E}(Z_{1}\ldots Z_{p})\leq \mathbb{E}(Z_{1})\ldots \mathbb{E}%
(Z_{p})+(p-1)\alpha \Vert Z_{1}\Vert _{\infty }\cdots \Vert Z_{p}\Vert
_{\infty }\,,
\end{equation*}%
and 
\begin{equation*}
\mathbb{E}(Z_{1})\ldots \mathbb{E}(Z_{p})\leq \mathbb{E}(Z_{1}\ldots
Z_{p})+(p-1)\alpha \Vert Z_{1}\Vert _{\infty }\cdots \Vert Z_{p}\Vert
_{\infty }\,.
\end{equation*}
\end{lemma}

\noindent {\bf Acknowledgments.} The authors are indebted to the referee for carefully reading of the manuscript and for helpful comments.

\end{document}